\documentclass[12pt]{amsart}
\usepackage{amscd,amssymb}
\usepackage[arrow,matrix]{xy}
\usepackage[plainpages,backref,urlcolor=blue]{hyperref}

\theoremstyle{plain}
\newtheorem{thm}[subsection]{Theorem}
\newtheorem{lem}[subsection]{Lemma}

\newtheorem{cor}[subsection]{Corollary}

\theoremstyle{definition}
\newtheorem{rk}[subsection]{Remark}
\newtheorem{definition}[subsection]{Definition}
\newtheorem{ex}[subsection]{Example}

\numberwithin{equation}{section}
\newcommand{\OO}{{\mathcal O}}

\newcommand{\I}{{\mathcal I}}

\newcommand{\al}{{\alpha}}

\newcommand{\Z}{\mathbb{Z}}

\newcommand{\C}{\mathbb{C}}

\newcommand{\PP}{\mathbb{P}}

\newcommand{\N}{\mathbb{N}}

\DeclareMathOperator{\Hom}{Hom}

\DeclareMathOperator{\im}{im}
\DeclareMathOperator{\coker}{coker}

\DeclareMathOperator{\defect}{def}

\begin{document}
%\date{June 4, 2009}

\title [Jacobian syzygies, stable sheaves and Torelli properties]
{Jacobian syzygies, stable reflexive sheaves, and Torelli properties for projective hypersurfaces with isolated singularities}

\author[Alexandru Dimca]{Alexandru Dimca$^1$}
\address{Univ. Nice Sophia Antipolis, CNRS,  LJAD, UMR 7351, 06100 Nice, France. }
\email{dimca@unice.fr}

\thanks{$^1$ Partially supported by Institut Universitaire de France}



\begin{abstract} We investigate the relations between the syzygies of the Jacobian ideal  of the defining equation for a projective hypersurface $V$ with isolated singularities  and the Torelli properties  of $V$ (in the sense of Dolgachev-Kapranov). We show in particular that hypersurfaces with a small Tjurina numbers are Torelli in this sense. When $V$ is a plane curve, or more interestingly, a surface in $\PP^3$, we  discuss the stability of the reflexive sheaf of logarithmic vector fields along $V$. A new lower bound for the minimal degree of a syzygy associated to a 1-dimensional complete intersection is also given.

\end{abstract}

\maketitle

%\tableofcontents

\section{Introduction} \label{sec:intro}

Let $X$ be the complex projective space $\PP^n$ and consider the associated graded $\C$-algebra $S= \oplus _kS_k $, with $S_k=H^0(X, \OO_X(k))$. For a nonzero section $f \in S_N$ with $N>1$, thought of as a homogeneous polynomial of degree $N$, we consider the hypersurface $V=V(f)$ in $X$ given by the zero locus of $f$and let $Y$ denote the singular locus of $V$, endowed with its natural scheme structure, see \cite{DBull}.
We assume in this paper that $V$ has isolated singularities.

Let $\I_Y \subset \OO_X$ be the ideal sheaf defining this 0-dimensional subscheme $Y \subset X$ and consider the graded ideal $I= \oplus_k I_k$ in $S$ with $I_k=H^0(X, \I_Y(k))$. Let $Z=Spec(S)$ be the corresponding affine space $\C^{n+1}$ and denote by $\Omega^k=H^0(Z, \Omega_Z^k)$ the $S$-module of global, regular $k$-forms on $Z$. 
Using a linear coordinate system $x=(x_0,...,x_n)$ on $X$, one sees that there is a natural
grading on $\Omega^k$, see \cite{DSt14} for details if necessary.

There is a well defined differential $1$-form $df \in \Omega^1$ and using it  we define two graded $S$-submodules in $\Omega^n$, namely
\begin{equation} 
\label{eqAR1}
AR(f)=\ker \{df \wedge : \Omega^n \to \Omega^{n+1}\}
\end{equation}
and
\begin{equation} 
\label{eqKR}
KR(f)=\im \{df \wedge : \Omega^{n-1}\to   \Omega^n\}.
\end{equation}
If one computes in a coordinate system $x$, then $AR(f)_m$ is the vector space of {\it all relations} of the type 
\begin{equation} 
\label{eqR}
R_m: a_0f_{x_0}+...a_nf_{x_n}=0,
\end{equation}
 with $f_{x_j}$ being the partial derivative of the polynomial $f$ with respect to $x_j$ and $a_j \in S_m$. Moreover, $KR(f)$ is the module of {\it Koszul relations} spanned by obvious relations of the type $f_{x_j}f_{x_i}+(-f_{x_i})f_{x_j}=0$ and the quotient
\begin{equation} 
\label{eqER}
ER(f)=AR(f)/KR(f)
\end{equation}
is the graded module of {\it essential relations} (which is of course nothing else but the $n$-th cohomology group of the Koszul complex of $f_{x_0},...,f_{x_n}$, maybe up to a shift in grading), see \cite{DBull}, \cite{DSt}. Note also that with this notation, the ideal $I$ is just the saturation of the Jacobian ideal $J_f=(f_{x_0},...,f_{x_n})\subset S=\C[x_0,...,x_n]$.

Let $\alpha_V$ be the Arnold exponent of the hypersurface $V$, which is by definition the minimum of the Arnold exponents of the singular points of $V$, cite \cite{DS1}, \cite{DS2}. Using Hodge theory, one can prove that
\begin{equation} 
\label{van1}
   ER(f)_m= 0 \text{ for  any } m <\alpha_VN -n,
\end{equation}
under the {\it additional hypothesis that all the singularities of $V$ are weighted homogeneous}, see \cite{DS2} and \cite{Kl}. It is interesting to note that even though the approaches in \cite{DS2} and \cite{Kl} are quite different, the condition that the singularities of $V$ are weighted homogeneous plays a key role in both papers.

While this inequality is the best possible in general, as one can see by considering hypersurfaces with a lot of singularities, see \cite{DiSt2}, \cite{Di4}, for  situations when the hypersurface $V$ has a small number of singularities this result is far from optimal.
Our first result gives the following  better bound in this case.

\begin{thm}
\label{thm1I}
Assume that the hypersurface $V:f=0$ in $\PP^n$ has degree $N$ and only isolated singularities. Then
$$ER(f)_m=0 \text{ for  any } m \leq n(N-2)-\tau(V),$$
where $\tau(V)$, the Tjurina number of $V$, is the sum of the Tjurina numbers of all the singularities of $V$.
\end{thm}
See also Theorem \ref{thm1} for a stronger result. The proof of these results is elementary (i.e. without Hodge theory) and it does not require the hypothesis of $V$ having weighted homogeneous singularities. In fact it applies to the following more general situation.

Let ${\bf f}=(f_0,...,f_n)$ be a collection of $n+1$ homogeneous polynomials in $S_e$ for some degree $e>0$. Assume that the associated ideal $J_{\bf f}=(f_0,...,f_n)$ in $S$ is a 
{\it 1-dimensional complete intersection}, i.e. $\dim S/J_{\bf f}=1$ and $(f_1,...,f_n)$ is a regular sequence in $S$. Replacing the differential $df$ by the 1-form 
$$\omega= f_0dx_0+...+f_ndx_n,$$
one can define the graded $S$-module $AR({\bf f})$ (resp. $ER({\bf f})$) by the formula \eqref{eqAR1} (resp. \eqref{eqER}). The ideal $J_{\bf f}$ defines a 0-dimensional subscheme $Y({\bf f})$ of $\PP^n$ and for each point $p$ in the support $|Y({\bf f})|$ of this scheme, we consider the corresponding local ring $\OO_{Y({\bf f}),p}$. This is an Artinian local ring, with maximal ideal denoted by $m_p$, and we denote by $o(m_p)$ the order of this ideal, i.e. the smallest integer $k>0$ such that $m_p^k=0$ in $\OO_{Y({\bf f}),p}$. With this notation we have the following
generalization of Theorem \ref{thm1}.

\begin{thm}
\label{thm1II}
Assume the ideal $J_{\bf f}=(f_0,...,f_n)$ in $S$ is a 
 1-dimensional complete intersection with $f_j \in S_e$ for all $j=0,...,n.$   Then
$$ER({\bf f})_m=0 \text{ for  any } m \leq n(e-1)-\sum _ {p\in |Y({\bf f})|} o(m_p).$$
\end{thm}
Let us come back to the case of a projective hypersurface with isolated singularities. The exact sequence of coherent sheaves on $X$ given by
\begin{equation} 
\label{es1}
 0 \to T\langle V\rangle \to \OO_X(1)^{n+1} \to \I_Y(N) \to 0,
\end{equation}
where the last non-zero morphism is induced by $(a_0,...,a_n) \mapsto a_0f_{x_0}+...a_nf_{x_n}$
can be used to define the sheaf $T\langle V\rangle$ of logarithmic vector fields along $V$, see
\cite{Se}. This is a reflexive sheaf, in particular a locally free sheaf $T\langle V\rangle$ (identified to a rank two vector bundle on $X$) in the case $n=2$. The above exact sequence clearly yields
\begin{equation} 
\label{eqAR}
AR(f)_m=H^0(X,T\langle V\rangle(m-1)),
\end{equation} 
for any integer $m$. This equality can be used to show the the reflexive sheaf $T\langle V\rangle$ is stable in many cases. This was done already in the case $n=2$ in \cite{DS14} and the corresponding result is stated below in Corollary \ref{cor2} and Example \ref{T-curves} without the hypothesis of weighted homogeneous singularities. But now we can go to higher dimension, and as an illustration we prove the following result.

\begin{thm}
\label{thm3}
Assume that the surface $V:f=0$ in $\PP^3$ has degree $N=3m+4\geq 4$ and only isolated singularities such that $$\tau(V) \leq 8m+5.$$
Then $F=T\langle V\rangle(m)$ is a normalized stable rank 3 reflexive sheaf on $\PP^3$ with Chern classes
$$c_1(F)=0, \  \ c_2(E)=6(m+1)^2 \text{  and  } c_3(E)=-23m^3-60m^2-60m-20+\tau(V).$$
\end{thm}
Similar results can be proved for $N$ congruent to 0 or 2 modulo 3, but the details are left to the interested readers. Recall the following notion.

\begin{definition}\label{deftor}
A reduced hypersurface $V\subset X=\PP^n$ is called  \emph{DK-Torelli} (where DK stands for Dolgachev-Kapranov) if   the hypersurface $V$ can be reconstructed as a subset of $X$ from the sheaf
$T\langle V \rangle$.
 \end{definition}
For a discussion of this notion and various examples we refer to \cite{DK}, \cite{UY}, \cite{DS14}.
In particular,  E. Sernesi and the author have shown in \cite{DS14} that the nodal curves with a small number of nodes are DK-Torelli. In the proof, which follows the line of the proof for smooth hypersurfaces outlined by K. Ueda and  M. Yoshinaga in \cite{UY}, we have used the inequality \eqref{van1}
for $n=2$. Since in the case of small number of singularities the bound obtained in Theorem \ref{thm1I} is better, it is natural to see if this new bound gives a slithly stronger result.
 Theorem \ref{thm1I} is hence applied to  prove the following result, which slightly improves in  the results on the Torelli properties of nodal (or nodal and cuspidal) curves obtained in a recent joint work with E. Sernesi, see \cite{DS14}. The following result also extends the result by K. Ueda and  M. Yoshinaga concerning smooth hypersurfaces in \cite{UY} to hypersurfaces having a small Tjurina number.

\begin{thm}
\label{thm2}
Let $V:f=0$ be a degree $N \geq 4$ hypersurface in $\PP^n$, having  only isolated singularities.
 If 
$$\tau(V) \leq \frac{(n-1)(N-4)}{2}+1,$$
 then one of the following holds.
\begin{enumerate}

\item $V$ is DK-Torelli;

\item $V$ is of Sebastiani-Thom type, i.e. in some linear coordinate system $(x_0,...,x_n)$ on $\PP^n$, the defining polynomial $f$ for $V$ is written as a sum $f=g+h$, with $g$ (resp. $h$) a polynomial involving only $x_0,...,x_r$ (resp. $x_{r+1},...,x_n$) for some integer $r$ satisfying $0 \leq r<n$.

\end{enumerate}

\end{thm}

The interest in having both Theorem \ref{thm3} and Theorem \ref{thm2} is that this allows the construction of  injective mappings from the varieties parametrizing surfaces with a fixed type of singularities into the moduli spaces of rank 3 stable reflexive sheaves on $\PP^3$, exactly as in the case $n=2$ discussed in \cite{DS14}. The case of rank 2 stable bundles on surfaces is rather well understood, see for instance \cite{HL},  but the case of rank 3 reflexive sheaves on $\PP^3$ seems to be still very mysterious.
Our construction may bring new light in the study of such moduli spaces, a problem which have attracted a lot of attention in the past, see \cite{Coanda}, \cite{Har}, \cite{Maru}, \cite{MR1}, \cite{Ok}, \cite{OSS}.

\section{The new bound on the minimal degree of a  syzygy}

Let $\OO_n$ denote the ring of holomorphic function germs at the origin of $\C^n$ and let 
$ m_n \subset \OO_n$ be its unique maximal ideal. For a function germ $g \in \OO_n$ defining an isolated hypersurface singularity at the origin of $\C^n$, we introduce an invariant
\begin{equation} 
\label{order}
 a(g)=\min \{ a \in \N \\ : \\m_n^a \subset   J_g+(g) \},
\end{equation}
where $J_g$ is the Jacobian ideal of $g$ in $\OO_n$ and $(g)$ is the principal ideal spanned by $g$ in $\OO_n$. We let $M(g)=\OO_n/J_g$ denote the {\it Milnor algebra} of the singularity $g$, and
$T(g)=\OO_n/(J_g+(g))$ the corresponding {\it Tjurina algebra}. Note that $a(g)=o(\tilde m_0)$, where $\tilde m_0$ denotes the maximal ideal of the Tjurina algebra $T(g)$. This invariant $a(g)$ clearly depends only on the contact class of the germ $g$ (i.e. the isomorphism class of the analytic germ $(\{g=0\},0)$), and this is the reason why we keep this notation besides the notation $o(\tilde m_0)$. This invariance is crucial in the computations given in the following example.

\begin{ex}
\label{exorder} (i) If $g=0$ is a node, i.e. an $A_1$-singularity, then $a(g)=1$.\\
(ii) If $g=0$ is a cusp, i.e. an $A_2$-singularity, then $a(g)=2$.\\
(iii) If $g=0$ is a $D_4$-singularity, e.g. an ordinary 3-tuple point when $n=2$, then $a(g)=3$.\\
(iv) One always has $a(g) \leq \tau (g)$, where $\tau (g)= \dim T(g)$ is the Tjurina number of $g$. Usually this inequality is strict, for instance when $g=x^d+y^d$ is an ordinary point of multiplicity $d$ and $n=2$, one has $a(g)=2d-3 <(d-1)^2$ for $d \geq 3$. The case $d=3$ corresponds to the $D_4$-singularity in dimension 2 mentioned above.

\end{ex}
One has a natural morphism $\OO_X(k) \to \OO_X(k)/\I_Y(k)$ for any integer $k$, inducing an evaluation morphism
\begin{equation} 
\label{eval}
 ev_k: S_k=H^0(X, \OO_X(k)) \to H^0(X, \OO_X(k)/\I_Y(k))=H^0(Y, \OO_Y),
\end{equation}
where the last equality comes from the fact that $Y$ is $0$-dimensional, i.e. its support $|Y|=V_{sing}$, the singular set of $V$, consists of finitely many points $p_1,...,p_s$. This fact also implies
\begin{equation} 
\label{eval2}
 H^0(Y, \OO_Y)=\oplus_{p \in V_{sing}}\OO_{Y,p}.
\end{equation}
On the other hand, if $g_p=0$ is a local analytic equation for the hypersurface singularity $(V,p)$, one has an isomorphism $\OO_{Y,p}=T(g_p)$ of local $\C$-algebras. The following result is elementary and it has appeared in various forms, see for instance Corollary 2.1 in \cite{BeSo}, Proposition (1.3.9) in \cite{D0}, or section 3 in \cite{Me}.
\begin{lem}
\label{vanish}
The evaluation morphism
$ev_k:S_k \to \oplus_{p \in V_{sing}}\OO_{Y,p}$
is surjective for any $k \geq \sum_{p \in V_{sing}} a(g_p)-1$. In other words, if one defines the $k$-th defect of the singular locus subscheme $Y$ by
$$\defect_kY=\dim \coker ev_k,$$
then $\defect_kY=0$ for $k \geq \sum_{p \in V_{sing}} a(g_p)-1$.
\end{lem}

\proof One considers the following decomposion of the evaluation map $ev_k$
$$S_k \to \oplus_{p \in |Y|}\OO_{X,p}/m_{n,p}^{a(g_p)} \to \oplus_{p \in |Y|}\OO_{Y,p},$$
with $m_{n,p}$ the maximal ideal of  $\OO_{X,p}$.  Then one notices that the first morphism is surjective by Corollary 2.1 in \cite{BeSo}, and the second morphism is surjective by the definition of the invariants $a(g_p)$.

\endproof

When $J_{\bf f}=(f_0,...,f_n)$ in $S$ is a 
 1-dimensional complete intersection as in the Introduction, we have the following similar result
(with an identical proof).

\begin{lem}
\label{vanish2}
The evaluation morphism
$ev_k:S_k \to \oplus_{p\in |Y({\bf f})|}\OO_{ Y({\bf f}),p} $
is surjective for any $k \geq \sum_{p\in |Y({\bf f})|   } o(m_p)-1$. In other words, if one defines the $k$-th defect of the  subscheme $ Y({\bf f}) $ by
$$\defect_kY({\bf f})=\dim \coker ev_k,$$
then $\defect_kY({\bf f})=0$ for $k \geq \sum_{p\in |Y({\bf f})|   } o(m_p)-1$.
\end{lem}

The main result of this section in the hypersurface case is the following.
\begin{thm}
\label{thm1}
Assume that the hypersurface $V:f=0$ in $\PP^n$ has degree $N$ and only isolated singularities, with local equations $g_p=0$ for $p \in V_{sing}.$ Then
$$ER(f)_m=0 \text{ for  any } m \leq n(N-2)-\sum_{p \in V_{sing}} a(g_p).$$
\end{thm}

\proof Using Theorem 1 in \cite{DBull}, we see that
\begin{equation} 
\label{relvsdef}
\dim ER(f)_m= \defect_{nN-2n-1-m}Y.
\end{equation}
The claim follows then from Lemma \ref{vanish}.
\endproof

Recall that for a homogeneous polynomial $f \in S$ we define its Milnor (or Jacobian) graded algebra to be the quotient $M(f)=S/J_f$. Then the {\it coincidence threshold} $ct(V)$ was defined as
$$ct(V)=\max \{q~~:~~\dim_K M(f)_k=\dim M(f_s)_k \text{ for all } k \leq q\},$$
with $f_s$  a homogeneous polynomial in $S$ of degree $N$ such that $V_s:f_s=0$ is a smooth hypersurface in $\PP^n$. Finally, the {\it minimal degree of a nontrivial relation} $mdr(V)$ is defined as
$$mdr(V)=\min \{q~~:~~ ER({ f})_{q}\ne 0\}.$$
It is known that one has the equality
\begin{equation} 
\label{REL}
ct(V)=mdr(V)+N-2,
\end{equation} 
see \cite{DSt}, formula (1.3). Theorem \ref{thm1} clearly implies the following.

\begin{cor}
\label{cor1}
Assume that the hypersurface $V:f=0$ in $\PP^n$ has degree $N$ and only isolated singularities, with local equations $g_p=0$ for $p \in V_{sing}.$ Then
$$mdr(V)  \geq n(N-2)-\sum_{p \in V_{sing}} a(g_p)+1$$
and
$$ct(V)  \geq T-\sum_{p \in V_{sing}} a(g_p)+1,$$
with $T=(n+1)(N-2)$.
\end{cor}

\begin{ex}
\label{nodal} Consider a nodal hypersurface $V$ in $\PP^n$ having $\sharp A_1$ singularities $A_1$. In this case $\alpha_V=n/2$, see  \cite{DS2}, hence the inequality \ref{van1} yields
$$ER(f)_m=0 \text{ for  any } m <n(N-2)/2.$$
On the other hand, Theorem \ref{thm1} yields
$$ER(f)_m=0 \text{ for  any } m \leq n(N-2)-\sharp A_1.$$

The second vanishing result is stronger than the first one exactly when
$n(N-2)/2 \leq n(N-2)-\sharp A_1$, i.e. if and only if 
$$\sharp A_1 \leq n(N-2)/2.$$
For $\sharp A_1=1$, this implies $ct(V)\leq n(N-2)+N-2=(n+1)(N-2)=T$ and we know that this is in fact an equality by Example 4.3 (i) in \cite{DSt}. Similarly, Example 4.3 (ii) in \cite{DSt} shows that
$ct(V)=T-1$ when $\sharp A_1=2$. Hence in these two cases the inequality in Theorem 
\ref{thm1} is in fact an equality. Example 4.3 (iii) in \cite{DSt} shows that
$ct(V)=T-1$ or $ct(V)=T-2$ when $\sharp A_1=3$, depending on whether the three nodes are collinear or not. It follows that the bound given by Theorem \ref{thm1} is optimal for $\sharp A_1 \leq 3$.
\end{ex}

\begin{ex}
\label{curves} Consider a reduced plane curve $V:f=0$ in $\PP^2$ having $n_k$ ordinary singularities  of multiplicity $k$ for $k=2,3,4$ and no other singularities. 
Theorem \ref{thm1} and Example  \ref{exorder} yield
$$ER(f)_m=0 \text{ for  any } m \leq 2(N-2)-n_2-3n_3-5n_4.$$
In the nodal case, i.e. when $n_3=n_4=0$, this bound can be better than the one given by 
the inequality \ref{van1}, but only when $V$ is irreducible (indeed, otherwise $ER(f)_{N-2} \ne 0$ as shown in \cite{DSt} via Hodge theory and in \cite{EU} without Hodge theory and in a more general setting).
\end{ex}

Now we consider the case of  a 1-dimensional complete intersection and give {\bf the proof of Theorem}
\ref{thm1II}. The proof is very similar to that of  Theorem
\ref{thm1I}, but we choose to present these results separatedly in order to preserve the additional geometric flavor of  Theorem
\ref{thm1I}. Note the the key equation \eqref{relvsdef} holds in the more general case of a 1-dimensional complete intersection as noted at the end of Remark 4 in \cite{DBull}. Namely, we have
\begin{equation} 
\label{relvsdef2}
\dim ER({\bf f})_m= \defect_{n(e-1)-1-m}Y({\bf f}).
\end{equation}

It remains to apply Lemma \ref{vanish2} instead of Lemma \ref{vanish} to complete the proof of Theorem
\ref{thm1II}. Here is one easy consequence.

\begin{cor}
\label{cor1.5}
Let ${\bf g}=(g_1,...,g_n)$ be a collection of $n$ homogeneous polynomials in $S_e$ for some degree $e>0$ such that $(g_1,...,g_n)$ is a regular sequence in $S$.  The ideal spanned by the $g_i$'s define a 0-dimensional subscheme $Y({\bf g})$ and for each point $q$ in the support 
$|Y({\bf g})|$ of this scheme, let $m'_q$ be the maximal ideal in the local ring 
$\OO_{Y({\bf g}),q}$ and let $o(m'_q)$ denote its order. 
Then $$\sum_{q \in |Y({\bf g})|}o(m'_q) \geq n(e-1)+1.$$
\end{cor} 

\proof Apply Theorem \ref{thm1II} to the collection ${\bf f}=(f_0,...,f_n)$, with $f_0=0$ and $f_i=g_i$ for $i>0$. The relation $1 \cdot f_0=0$ implies that $ER({\bf f})_0 \ne 0$, and this yields the result.

\endproof

\begin{rk}
\label{rk1}
(i) The example $g_i=x_i^ e$ for $i=1,...,n$, shows that the inequality in Corollary \ref{cor1.5} is sharp. 

(ii) It is  clear that $o(m'_q) \leq \dim_{\C } \OO_{Y({\bf g}),q}$, where the equality holds if and only if the 0-dimensional singularity $(Y({\bf g}),q)$ has embedding dimension at most one.
In particular we get
$$\sum_{q \in |Y({\bf g})|}o(m'_q) \leq \sum_{q \in |Y({\bf g})|}\dim_{\C } \OO_{Y({\bf g}),q}=e^n.$$

(iii) The reader interested in the case of  a 1-dimensional complete intersection 
${\bf f}=(f_0,...,f_n)$, where the degrees of the $f_i$'s are different, can also obtain bounds on the degree of syzygies following the above approach. Indeed, the relation \eqref{relvsdef2} essentially continues to hold as noted at the end of Remark 4 in \cite{DBull} and the references given there, but additional care is needed to get the right homogeneous components in $ER({\bf f})$ in this case.
\end{rk}

\section{Stability of $T\langle V\rangle$  for $n=2$ and $n=3$} \label{sec4}
First we discuss the simpler case $n=2$.
Using Proposition 2.4 in \cite{Se} which says that  $T\langle V\rangle$ is stable if and only if
$AR(f)_m=0$ for all $m \leq   (N-1)/2$,
we get the following consequence of our Theorem \ref{thm1}.

\begin{cor}
\label{cor2}
Assume that the curve $V:f=0$ in $X=\PP^2$ has degree $N$ and only isolated singularities, with local equations $g_p=0$ for $p \in V_{sing}$. Then the vector bundle $T\langle V\rangle$ is stable if 
\begin{equation} 
\label{eqstab}
\left [ (N-1)/2 \right ]  \leq 2(N-2)-\sum_{p \in V_{sing}} a(g_p),
\end{equation} 
where $[y]$ denotes the largest integer verifying $[y]\leq y$.
In particular, since a rank two stable vector bundle is not splittable, it follows that $V$ is not a free divisor when the  inequality \eqref{eqstab} holds.
\end{cor}

\begin{ex}
\label{T-curves} For $N \geq 5$, consider the family of  plane curves $V_N:f_N=0$ in $\PP^2$ given by the equation
$$f_N=x^2y^2z^{N-4}+x^5z^{N-5}+y^5z^{N-5}+x^N+y^N=0.$$
Then $V_N$ has a unique singularity at $p_1=(0,0,1)$ which is isomorphic to the singularity
$g(u,v)=u^2v^2+u^5+v^5$. It follows that $\tau(V_N)=\tau(g)=10$ and $a(g)=5$, see for instance Example (6.56) in \cite{D00}. Moreover, this singularity , usually denoted by $T_{2,5,5}$ in Arnold's classification, is not weighted homogeneous, since $11=\mu(g)>\tau(g)=10.$
 
Theorem \ref{thm1}  yields $mdr(V_N)\geq 2N-8$,
while Theorem \ref{thm1I}  yields the weaker bound $mdr(V_N)\geq 2N-13.$
A direct computation of the Jacobian syzygies in the case $5 \leq N \leq 10$ using Singular shows that $mdr(V_N)=2N-7$. Therefore Theorem \ref{thm1} is almost sharp in these cases.

Using Corollary \ref{cor2}, this computation also implies that the curves $V_N$ have the property that the associated bundle $T\langle V_N\rangle$ is stable for any $N \geq 5$.
\end{ex}

\begin{ex}
\label{free} If $V$ is an irreducible free divisor in $X=\PP^2$ with degree $N$ and only isolated singularities, with local equations $g_p=0$ for $p \in V_{sing}$, it follows from Corollary  \ref{cor2}
that one has 
\begin{equation} 
\label{eqfree}
\sum_{p \in V_{sing}} a(g_p) > 2(N-2)-\left [ (N-1)/2 \right ]
\end{equation} 
In other words, such a curve should have a lot of singularities (or singularities with large invariants $a(g_i)$) and this explains the difficulty and the interest in constructing such examples, see 
for instance \cite{BC}, \cite{Na}, \cite{ST}. The example $V_5$ above shows that the inequality \eqref{free} in not sufficient to imply the freeness of the divisor.
\end{ex}

Now we pass to the case $n=3$ and prove Theorem \ref{thm3}. First we compute the Chern classes of the sheaf 
$T\langle V\rangle$ using the exact sequence \eqref{es1} which yields
$$c(T\langle V\rangle)\cdot c(\I_Y(N))=(1+\al)^4$$
where $\al \in H^2(\PP^3)$ is the standard generator and $\al^4=0$.
To compute $c(\I_Y(N))$ we use the exact sequence
\begin{equation} 
\label{esY}
 0 \to \I_Y(N) \to \OO_X(N) \to \OO_{Y} \to 0
\end{equation} 
which gives
$$c(\I_Y(N))\cdot (1+\tau(V)\al^3)=1+N\al.$$
Using this, we finally get by a direct computation the following Chern classes (identified to the corresponding Chern numbers)
$c_1(T\langle V\rangle)=-N+4$, $c_2(T\langle V\rangle)=N^2-4N+6$ and
$$c_3(T\langle V\rangle)=-N^3+4N^2-6N+4+\tau(V).$$
To study the stability of $T\langle V\rangle$, we have first to normalize it, i.e. find an integer $m$ such that $c_1(T\langle V\rangle(m)) \in \{0,-1,-2\}$.
Now $c_1(T\langle V\rangle(m))=-N+4+3m$ and we discuss below only the case when $N$ is congruent to $1$ modulo 3, i.e. there is a unique $m$ such that $N=3m+4$. However, everything that follows works identically for the other two cases. Then using Remark 1.2.6 (b) in \cite{OSS}, we have that $F=T\langle V\rangle(m)$ is stable if and only if one has
\begin{equation} 
\label{stab}
H^0(X, F)=H^0(X,F^*)=0,
\end{equation} 
where $X=\PP^3$ and $F^*$ is the dual of $F$. The first condition is easy to check. Indeed, using \eqref{eqAR}, we have
$H^0(X,F)=AR(f)_{m+1}$. On the other hand, $AR(f)_{m+1}=0$ since the condition
$\tau(V) \leq 8m+5$ is exactly what we need to apply Theorem \ref{thm1I}. Note also that
$m+1 =(d-1)/3<d-1$ and hence $AR(f)_{m+1}=ER(f)_{m+1}$.
The second vanishing $H^0(X,F^*)=0$ requires more work.
If we dualize the exact sequence \eqref{es1}, we get
\begin{equation} 
\label{es2}
0 \to \I_Y(N)^* \to \OO_X(-1)^4 \to T\langle V\rangle^* \to {\mathcal Ext}^1( \I_Y(N),\OO_X) \to...
\end{equation} 
We  show that ${\mathcal Ext}^1( \I_Y(N),\OO_X)=0$ by proving the corresponding vanishing at the stalk level. There are three cases to discuss.

\noindent{\bf  Case 1.} If $x \in X$ is not a singular point in $Y$, then $  \I_Y(N)_x=\OO_{X,x}$ and hence clearly $Ext^1(\I_Y(N)_x,\OO_{X,x})=0$.

\noindent {\bf Case 2.}  Assume $x \in X$ is  a weighted homogeneous singular point of $V$, with local equation 
$g=0$. Then $g \in J_g$ and hence $  \I_Y(N)_x=J_g$, which is a 0-dimensional complete intersection. The exact sequence
$$0 \to J_g \to \OO_n \to M(g) \to 0$$
yields a long exact sequence containing the sequence
$$0=Ext^1(\OO_n,\OO_n) \to Ext^1(J_g, \OO_n) \to Ext^2(M(g),\OO_n)=0.$$
The last vanishing is a well known property of complete intersections, see for instance \cite{GH}, p. 690.

\noindent {\bf Case 3.}  Consider now the case of a singularity $x$ given by $g=0$ which is not weighted homogeneous and hence $  \I_Y(N)_x=(g)+J_g$. The exact sequence
$$ 0 \to J_g \to (g)+J_g \to ((g)+J_g)/J_g \to 0$$
and Case 2. show that it is enough to prove that $Ext^1(((g)+J_g)/J_g, \OO_n)=0.$
Let $K_g$ be the kernel of the morphism given by multiplication $g:M(g) \to M(g)$. Then we have a short exact sequence
$$ 0 \to K_g \to M(g) \to ((g)+J_g)/J_g \to 0$$
which gives the result since
$$0=\Hom (K_g, \OO_n) \to Ext^1(((g)+J_g)/J_g, \OO_n) \to Ext^1(M(g),\OO_n)=0.$$
In conclusion,  \eqref{es2} is a short exact sequence. Twisting by $(-m)$ we get
$$0 \to  \I_Y(N+m)^* \to \OO^4_X(-1-m) \to F^* \to 0.$$
To show that $H^0(X,F^*)=0$ it is enough to show that $H^1( \I_Y(N+m)^*)=0$.
Notice that
$$\I_Y(N+m)^*={\mathcal Hom}( \I_Y(N+m), \OO_X)={\mathcal Ext^0}( \I_Y(N+m-4), \omega_X)$$
where $\omega_X=\OO_X(-4)$ is the dualizing sheaf of $X$. Now use the spectral sequence
$$E_2^{p,q}=H^p(X,{\mathcal Ext^q}(\I_Y(N+m-4),\omega_X)$$
converging to $Ext^{p+q}(\I_Y(N+m-4),\omega_X)=H^{3-p-q}(X,\I_Y(N+m-4))$.
Twisting the exact sequence \eqref{esY} by $(m-4)$ we get $H^{2}(X,\I_Y(N+m-4))=0$ which implies via the spectral sequence $H^1( \I_Y(N+m)^*)=0$.
To complete the proof of Theorem \ref{thm3} it suffices use the formulas for the Chen classes of a tensor product, see \cite{OSS}, p. 16, and get
$$c_1(T\langle V\rangle(m))=-N+4+3m=0$$
by the choice of $m$,
$$c_2(T\langle V\rangle(m))=c_2(T\langle V\rangle)+2m c_1(T\langle V\rangle)+3m^2$$
and
$$c_3(T\langle V\rangle(m))=c_3(T\langle V\rangle)+mc_2(T\langle V\rangle)+m^2c_1(T\langle V\rangle)+m^3.$$
Using the formulas given above for $c_i(T\langle V\rangle)$ and replacing $N=3m+4$ yields the claimed formulas.

\section{DK-Torelli type properties for singular hypersurfaces} \label{sec5}

Now we turn to the proof of Theorem \ref{thm2} stated in the Introduction.
This proof follows closely the proof of the corresponding result in \cite{DS14}. We repeat below the main steps, for the reader's convenience and also to point out the new facts necessary to treat the $n$-dimensional case.

\begin{lem}
\label{step1}
With the above notation and hypothesis, the sheaf $T\langle V \rangle$ determines the  vector subspace $J_{f,N-1}\subset S_{N-1}$.
\end{lem}

\bigskip

To prove this Lemma, let $E:g=0$ be a (possibly nonreduced) hypersurface in $X=\PP^n$ of degree $N-1$. For any $k \in \Z$, consider the exact sequence
$$ 0 \to  \OO_{X}(k-N+1) \to  \OO_{X}(k) \to  \OO_{E}(k) \to 0,$$
where the first morphism is induced by the multiplication by $g$.
Tensoring this sequence of  sheaves by $T\langle V\rangle$, we get a new short exact sequence
$$0 \to T\langle V\rangle(k-N+1) \to T\langle V\rangle(k) \to  T\langle V\rangle (k)\otimes\OO_E\to 0.$$
The injectivity of the first morphism comes from the fact that $T\langle V\rangle(k-N+1)$ is a torsion free sheaf.
The associated long exact sequence of cohomology groups looks like
$$0 \to H^0(T\langle V\rangle(k-N+1)) \to H^0(T\langle V\rangle(k)) \to 
H^0 (T\langle V\rangle (k)\otimes\OO_E) \to $$ $$ \to H^1(T\langle V\rangle(k-N+1)) \to H^1(T\langle V\rangle(k)) \to\cdots$$
Then, using the formula \eqref{eqAR}, we see that 
$$\delta_k= \dim  H^0(T\langle V\rangle(k))-\dim H^0(T\langle V\rangle(k-N+1))= 
\dim AR(f)_{k+1}-\dim AR(f)_{k-N+2}$$
 depends only on $f$ but not on $g$. Next note that the morphism
$$ H^1(T\langle C\rangle(k-N+1)) \to H^1(T\langle C\rangle(k)) $$
in the above exact sequence can be identified, using the formulas (5) and (9) in \cite{Se}  with the morphism
$$g^*_{k+1}: (I/J_f)_{k+1} \to (I/J_f)_{k+N}$$
induced by the multiplication by $g$ (we recall that $I$ is the saturation of the Jacobian ideal $J_f$ ). The above proves the following equality.
\begin{equation} \label{keyequality}
\dim H^0 (T\langle C\rangle (k)\otimes\OO_E) = \delta_k + \dim \ker g^*_{k+1}.
\end{equation}
Let $m$ be the largest integer such that $2m \leq N-2$.
Since clearly $m<N-1$, it follows that $J_{f,m}=0$ and hence $g^*_m$ is defined on $I_m$. If $g \in J_f$, then clearly $g^*_m=0$, and hence its kernel has maximal possible dimension.

To complete the proof of Lemma \ref{step1} it is enough to show that the converse also holds.
To do this, we  show first that there are two elements $h_1,h_2 \in I_m$ having no irreducible factor in common. Otherwise, all the elements in $I_m$ are divisible by a homogeneous polynomial, and hence in particular one has $\dim I_m \leq \dim S_{m-1}$ which implies
\begin{equation} \label{keyineq}
\tau(V) \geq \dim S_m/I_m \geq {m+n \choose n} - {m+n-1 \choose n} =  {m+n-1 \choose n-1} . 
\end{equation}
One also has the inequality (perhaps well known)
\begin{equation} \label{easyineq}
 {m+n-1 \choose n-1} \geq (n-1)m+1, 
\end{equation}
which can be proved by looking at the subsets $E'$ of cardinal $n-1$ of a set $E$ which is a disjoint union $E=E_1 \cup E_2$, with $\sharp E_1=m$, $\sharp E_2=n-1$ and count how many subsets $E'$ satisfy $\sharp (E' \cap E_1)\leq 1$.
It follows that
$$2m \leq 2 \cdot \frac{\tau(V)-1}{n-1} \leq N-4,$$
a contradiction with the choice of $m$. This shows  that there are two elements $h_1,h_2 \in I_m$ having no irreducible factor in common.

Then $g^*_m=0$ implies $gh_1=\sum_{j=0,n} a_jf_{x_j}$ and $gh_2=\sum_{j=0,n} b_jf_{x_j}$
 for some polynomials $a_j,b_j \in S_m$.
It follows that
$$\sum_{j=0,n}(a_jh_2-b_jh_1)f_{x_j}=0.$$
Since
$$\sum_{i=1,s}a(g_i) \leq \tau(V)=\sum_{i=1,s} \tau(g_i) \leq \frac{(n-1)(N-4)}{2}+1\leq (n-1)(N-2)$$
it follows that 
$$2m \leq N-2 \leq n(N-2)-\sum_{i=1,s}a(g_i).$$
Theorem \ref{thm1}  implies that the only syzygy of degree $2m$ is the trivial one, i.e.
$a_jh_2=b_jh_1$ for any $j$. These relations imply that the polynomials $a_j$'s are all divisible by $h_1$ in $S$, and hence $g\in J_f$.

It follows that $g \in J_{f,N-1}$ if and only if 
$$\dim H^0 (T\langle V\rangle (m-1)\otimes\OO_E) = \delta_{m-1} + \dim I_m,$$
i.e. the sheaf $T\langle V\rangle$ determines the homogeneous component $J_{f,N-1}$ of the Jacobian ideal $J_f$, and this completes the proof of Lemma \ref{step1}.

To finish the proof of Theorem \ref{thm2}, it is enough to use Theorem 1.1 in Zhenjian Wang paper \cite{W}, which generalizes a Lemma in \cite{DS14} covering the case $n=2$.
Indeed, this Theorem says that we can have the following siuations.

(A) The Jacobian ideal $J_f$ (or its homogeneous component $J_{f,N-1}$ determines $f$ up to a multiplicative nonzero constant. In this case $V$ is DK-Torelli.

(B) $V$ is of Sebastiani-Thom type.

(C)  $V$ has at least one singular point $p_i$ with multiplicity $N-1$, i.e. there is a local equation $g_i$ such that $g_i \in m_{n,p_i}^{N-1}$. However this is impossible in our conditions as we show now. If $g_i \in m_{n,p_i}^{N-1}$, then the monomials in the corresponding local coordinates $u_1,...,u_n$ of degree $\leq N-3$ are linearly independent in $\OO_n/(J_{g_i}+(g_i))$. This implies that
$$\tau(V) \geq \tau(g_i) \geq \dim S_{N-3}= {N-3+n \choose n} \geq n(N-3)+1,$$
as in \eqref{easyineq}. But this is in contradiction with the hypothesis
$$\tau(V) \leq \frac{(n-1)(N-4)}{2}+1,$$
so the proof of  Theorem \ref{thm2} is complete.

\begin{cor}
\label{cor3}
Assume that the  curve $V:f=0$ in $\PP^2$ has degree $N\geq 4$, $\nu$ nodes, $\kappa$ cusps and no other singularities. If
$$\nu +2\kappa \leq \frac{N-2}{2}$$
then the curve $V$ is DK-Torelli.
\end{cor}
This is a direct consequence of Theorem \ref{thm2}, since a curve with only nodes and cusps and of degree at least 4 cannot satisfy the property (2), see also \cite{DS14}.

For $V$ irreducible and $\kappa=0$, this result coincides with the result given in \cite{DS14}.
For the remaining cases,  Corollary \ref{cor3} is a slight improvement over the corresponding results given in \cite{DS14}. In particular, Corollary \ref{cor3} shows that a curve with $\nu=0$ and $\kappa=1$
is DK-Torelli as soon as $N \geq 6$, while the bound given in \cite{DS14} for the same result was $N\geq 8$.

\end{document}